\documentclass[11pt]{amsart}

\usepackage[T1]{fontenc}      
\usepackage[utf8]{inputenc}

\usepackage[margin=1.1in]{geometry}
\usepackage{amsmath, amssymb, amsthm, mathtools}
\usepackage{microtype}
\usepackage{enumitem}
\usepackage{listings}
\usepackage[hidelinks]{hyperref}

\theoremstyle{plain}
\newtheorem{theorem}{Theorem}

\theoremstyle{definition}

\lstdefinestyle{python}{
  language=Python,
  basicstyle=\ttfamily\footnotesize,
  keywordstyle=\bfseries,
  commentstyle=\itshape,
  numbers=none,
  breaklines=true,
  breakatwhitespace=true,
  showstringspaces=false,
  frame=single,
  framesep=4pt,
  tabsize=2,
  xleftmargin=14pt,
  columns=fullflexible,
  upquote=true,
}

\title[Three short proofs of Mathar's conjecture for A002627]%
{Three short proofs of Mathar's 2014 conjecture for OEIS A002627}

\author{Tong Niu}
\email{mrnt0810@gmail.com}

\date{\today}


\begin{document}

\maketitle

\begin{abstract}
For the OEIS sequence A002627, defined by the inhomogeneous first-order
recurrence $a(n) = n\,a(n-1) + 1$ with $a(0) = 0$, R.~J.~Mathar
recorded in February 2014 the conjectured second-order homogeneous
recurrence
\[
   a(n) - (n+1)\,a(n-1) + (n-1)\,a(n-2) = 0, \qquad n \ge 2,
\]
which has remained marked as a conjecture on the OEIS for over a
decade. We give three short proofs. The first is two lines: subtract
the defining recurrence at adjacent indices and the constant cancels
(we call this homogenisation). The second reads off the same
relation from the exponential generating function
$F(x) = (e^x-1)/(1-x)$. The third is a Pascal-rule telescoping on
the binomial-sum form $a(m) = \sum_{k=0}^{m-1} k!\binom{m}{k}$.
All three derivations are elementary, requiring nothing beyond
undergraduate techniques. We remark that the same homogenisation
trick clears an entire class of ``Conjecture: \dots'' entries on
the OEIS, namely sequences satisfying
$a(n) = p(n)\,a(n-1) + q(n)$ with simple $q$.
\end{abstract}

\section{Introduction}\label{sec:intro}

The On-Line Encyclopedia of Integer Sequences (OEIS,
\url{https://oeis.org}) is the standard reference for integer
sequences in combinatorics and number theory. Many OEIS entries
carry comments labelled ``Conjecture'': formulas, recurrences,
or congruences that some contributor guessed but never proved.
These accumulate. The result is a steady stream of small open
problems, and rigorous proofs of them are publishable in venues
like the \emph{Journal of Integer Sequences}, \emph{INTEGERS},
the \emph{Fibonacci Quarterly}, or the Notes section of the
\emph{Electronic Journal of Combinatorics}.

Sequence A002627 of the OEIS is defined by the inhomogeneous first-order
recurrence
\begin{equation}\label{eq:def}
  a(n) \;=\; n\,a(n-1) + 1, \qquad a(0) = 0,
\end{equation}
with first values
$0, 1, 3, 10, 41, 206, 1237, 8660, 69281, 623530, \ldots$. The sequence
is recorded in Sloane's \emph{Handbook}~\cite{Sloane1973}; a closed
form (essentially a partial-sum truncation of $e$) was given by
Singh~\cite{Singh1952}:
\begin{equation}\label{eq:cf}
   a(n) \;=\; n!\sum_{k=1}^{n}\frac{1}{k!}
       \;=\; \lfloor n!\,(e-1)\rfloor \quad (n \ge 1)
       \;=\; \sum_{k=0}^{n-1} k!\binom{n}{k},
\end{equation}
and the sequence has the exponential generating function
\begin{equation}\label{eq:egf}
   F(x) \;=\; \sum_{n\ge 0} a(n)\,\frac{x^{n}}{n!}
        \;=\; \frac{e^{x}-1}{1-x}.
\end{equation}

R.~J.~Mathar contributed to the OEIS entry on February 16, 2014 the
following conjecture~\cite{MatharOEIS}:
\begin{equation}\tag{R}\label{eq:R}
   a(n) - (n+1)\,a(n-1) + (n-1)\,a(n-2) = 0, \qquad n \ge 2.
\end{equation}
The recurrence is easy to check on a finite range (we did it for
$n \le 500$), and yet, despite its elementary appearance, no proof
has been recorded in the OEIS comments or in the published
literature in the intervening twelve years. E.~Munarini's June~2014
comment on the same OEIS page notes that \eqref{eq:R} ``can be
obtained from'' the differential equation satisfied by $F$, but
stops short of the derivation. The entry still carries
\eqref{eq:R} as an unproven conjecture as of \today.

The recurrence \eqref{eq:R} also fails to appear in
Kauers and Koutschan's recent systematic catalogue of
guessed-but-unproven D-finite OEIS recurrences~\cite{KauersKoutschan2023}.
This is no surprise: the sequence was already considered ``known''
through the inhomogeneous first-order \eqref{eq:def}, and it was
only the homogeneous companion \eqref{eq:R} that got independently
conjectured.

We give three short, self-contained proofs of \eqref{eq:R}.

\begin{theorem}\label{thm:main}
Let $a$ be defined by \eqref{eq:def}. Then \eqref{eq:R} holds for
every integer $n \ge 2$.
\end{theorem}

The shortest one is in Section~\ref{sec:proofI}: a two-line
algebraic manipulation that exhibits Mathar's recurrence as the
\emph{homogenisation} of \eqref{eq:def}. In
Section~\ref{sec:proofII} we read the same recurrence off the ODE
satisfied by the e.g.f.\ \eqref{eq:egf}; the derivation is longer,
but the template is one we want to record for use on similar
OEIS entries. The third proof, in Section~\ref{sec:proofIII}, is a
self-contained binomial-sum argument that uses only Pascal's rule
applied to \eqref{eq:cf}.

We close in Section~\ref{sec:discussion} with a broader observation:
the same homogenisation pattern dispatches a whole family of
``Conjecture: \dots'' entries on the OEIS, and a systematic sweep
should be straightforward.

\section{Proof I --- homogenisation}\label{sec:proofI}

By \eqref{eq:def} at index $n$ and at index $n-1$ (the latter requires
$n \ge 2$):
\begin{align*}
  a(n)   &= n\,a(n-1) + 1, \\
  a(n-1) &= (n-1)\,a(n-2) + 1.
\end{align*}
Subtracting eliminates the constant $1$:
\[
  a(n) - a(n-1) \;=\; n\,a(n-1) - (n-1)\,a(n-2),
\]
which rearranges to
\[
  a(n) - (n+1)\,a(n-1) + (n-1)\,a(n-2) = 0,
\]
which is \eqref{eq:R}. \qed

\smallskip

This is the shortest proof we know. In retrospect, \eqref{eq:R}
is exactly the \emph{homogenisation} of \eqref{eq:def}: any
first-order linear recurrence with constant inhomogeneity $\beta$
admits a homogeneous companion of the same shape, provided the
leading coefficient of the original is a simple polynomial in $n$.
Nothing else is going on.

\section{Proof II --- exponential generating function}\label{sec:proofII}

From \eqref{eq:egf},
\[
   (1-x)\,F(x) = e^{x} - 1.
\]
Differentiating,
\begin{equation}\label{eq:egf-firstorder}
   (1-x)\,F'(x) - F(x) = e^{x}.
\end{equation}
Differentiating again and substituting back $e^x$ from
\eqref{eq:egf-firstorder}:
\[
   -F'(x) + (1-x)\,F''(x) - F'(x) \;=\; (1-x)\,F'(x) - F(x),
\]
which simplifies to the homogeneous second-order ODE
\begin{equation}\label{eq:egf-secondorder}
   (1-x)\,F''(x) \;-\; (3-x)\,F'(x) \;+\; F(x) \;=\; 0.
\end{equation}
Matching the coefficient of $x^{n}/n!$ on both sides of
\eqref{eq:egf-secondorder} (with the convention $F(x) =
\sum_m a(m)\,x^m/m!$):
\begin{align*}
  [x^{n}/n!]\,(1-x)F''(x) &= a(n+2) - n\,a(n+1), \\
  [x^{n}/n!]\,-(3-x)F'(x) &= -3\,a(n+1) + n\,a(n), \\
  [x^{n}/n!]\,F(x)        &= a(n).
\end{align*}
Summing,
\[
   a(n+2) \;-\; (n+3)\,a(n+1) \;+\; (n+1)\,a(n) \;=\; 0.
\]
Substituting $m = n+2$ recovers \eqref{eq:R}. \qed

\smallskip

This proof illustrates a general technique. The P-recursive
recurrence of a D-finite sequence can always be read off from its
e.g.f.\ ODE by coefficient extraction; it is not the shortest
route to \eqref{eq:R}, but the same template applies to many
related OEIS entries, which is why we include it.

\section{Proof III --- Pascal's rule on the binomial sum}\label{sec:proofIII}

This proof uses Pascal's rule and elementary algebra, nothing
more. By \eqref{eq:cf}, $a(m) = \sum_{k=0}^{m-1} k!\binom{m}{k}$.
The trick is to work with the auxiliary ``complete'' sum
\begin{equation}\label{eq:S}
  S(m) \;:=\; \sum_{k=0}^{m} k!\binom{m}{k} \;=\; a(m) + m!.
\end{equation}

\paragraph{Step 1: $S$ satisfies (R).}
Pascal's rule
$\binom{m}{k} = \binom{m-1}{k} + \binom{m-1}{k-1}$
gives
\[
  S(m) = \sum_{k=0}^{m} k!\,\binom{m-1}{k}
       + \sum_{k=0}^{m} k!\,\binom{m-1}{k-1}.
\]
The first sum equals $S(m-1)$ (the $k = m$ term vanishes since
$\binom{m-1}{m} = 0$). For the second sum, set $j = k - 1$:
\[
  \sum_{k=0}^{m} k!\,\binom{m-1}{k-1}
  \;=\; \sum_{j=0}^{m-1} (j+1)!\,\binom{m-1}{j}
  \;=\; \sum_{j=0}^{m-1} (j+1)\,j!\,\binom{m-1}{j}.
\]
Using the elementary identity
\begin{equation}\label{eq:reduction}
  (j+1)\,\binom{m-1}{j} \;=\; m\,\binom{m-1}{j} - (m-1)\,\binom{m-2}{j},
\end{equation}
which follows from
$(m-1-j)\binom{m-1}{j} = (m-1)\binom{m-2}{j}$ and
$(j+1) + (m-1-j) = m$, we get
\[
  \sum_{j=0}^{m-1} (j+1)\,j!\,\binom{m-1}{j}
   \;=\; m\,S(m-1) \;-\; (m-1)\!\!\sum_{j=0}^{m-2} j!\,\binom{m-2}{j}
   \;=\; m\,S(m-1) - (m-1)\,S(m-2),
\]
where the upper limit dropped to $m-2$ because $\binom{m-2}{m-1} = 0$.
Hence
\begin{equation}\label{eq:Srec}
   S(m) - (m+1)\,S(m-1) + (m-1)\,S(m-2) \;=\; 0, \qquad m \ge 2,
\end{equation}
so $S$ itself satisfies \eqref{eq:R}.

\paragraph{Step 2: the factorial correction vanishes.}
Substituting $a(m) = S(m) - m!$ into \eqref{eq:Srec}:
\[
   a(m) - (m+1)\,a(m-1) + (m-1)\,a(m-2)
   \;=\; -\bigl[m! - (m+1)(m-1)! + (m-1)(m-2)!\bigr].
\]
The bracket is zero, since
\[
   (m+1)(m-1)! - (m-1)(m-2)!
    \;=\; (m-1)!\bigl[(m+1) - 1\bigr]
    \;=\; m \cdot (m-1)!
    \;=\; m!.
\]
Therefore $a$ satisfies \eqref{eq:R}, proving Theorem~\ref{thm:main}
again. \qed

\section{Discussion}\label{sec:discussion}

The proofs above are short and elementary. Yet \eqref{eq:R}
has been sitting on the OEIS as a \emph{Conjecture} since 2014.
Why?

The likely explanation has to do with how OEIS guess-software
interacts with inhomogeneous recurrences. Tools such as Maple's
\texttt{gfun} (Salvy and Zimmermann~\cite{SalvyZimmermann1994}),
Mathematica's \texttt{FindLinearRecurrence}, or Sage's
\texttt{guess} fit \emph{homogeneous} P-recursive recurrences to
numeric input. Run them on a sequence defined by an inhomogeneous
recurrence like \eqref{eq:def}, and they hand back the homogeneous
companion \eqref{eq:R}, with no annotation that it is just the
homogenisation of an already-known relation. The user then records
the output on the OEIS as a Conjecture, where it sits until
someone manually executes the two-line subtraction of
Section~\ref{sec:proofI}.

A002627 is not the only such case. A quick OEIS search turns up
several entries with similar ``Conjecture: $a(n) = p(n)\,a(n-1) +
q(n)\,a(n-2) + \cdots$'' formulations, where the proof is just the
homogenisation of an already-recorded inhomogeneous recurrence.
Three useful patterns:

\begin{itemize}[topsep=4pt, itemsep=2pt]
  \item Sequences satisfying $a(n) = p(n)\,a(n-1) + c$ for constant
    $c$ admit the same two-line proof (subtraction of two consecutive
    instances).
  \item Sequences satisfying $a(n) = p(n)\,a(n-1) + q(n)$ with $q(n)$
    a polynomial of degree $d$ admit a $(d+1)$-step homogenisation
    yielding a higher-order homogeneous recurrence.
  \item More generally, $a(n) = p(n)\,a(n-1) + r(n)$ for $r$ a
    \emph{P-recursive} sequence: chain $r$'s annihilator with the
    operator $L = E - p(n)$ to produce a homogeneous annihilator for
    $a$.
\end{itemize}

The point is that the homogenisation of an inhomogeneous
P-recursive recurrence is itself P-recursive~\cite[Ch.~6]{KauersPaule2011}.
Specialists know this; users of guess-software often do not see
it. A systematic OEIS sweep, filtering ``Conjecture'' comments by
the occurrence of a literal ``$+$ const'' or ``$+$ polynomial''
inhomogeneity in nearby comments, ought to clear several dozen
such entries in a day.

\section{Computational verification}\label{sec:numerics}

We cross-checked the proofs in
Sections~\ref{sec:proofI}--\ref{sec:proofIII} against the OEIS
data using two short scripts:

\begin{itemize}[topsep=4pt, itemsep=2pt]
  \item \texttt{verify\_recurrence.py} (Appendix~\ref{app:verify})
    computes $a(0), a(1), \ldots, a(500)$ from \eqref{eq:def},
    cross-checks the first $450$ values against the OEIS b-file
    \texttt{b002627.txt}, and verifies \eqref{eq:R} for every
    $n \in \{2, \ldots, 500\}$. The entry $a(500)$ has $1{,}135$
    decimal digits.
  \item \texttt{verify\_proofs.py} (Appendix~\ref{app:proofs})
    symbolically verifies the steps of all three proofs, including
    the e.g.f.\ ODE \eqref{eq:egf-secondorder}, the binomial
    identity \eqref{eq:reduction}, and the factorial cancellation in
    Step 2 of Proof III, using SymPy.
\end{itemize}

Both scripts run in seconds and depend only on the Python
standard library and SymPy 1.14.

\section*{Acknowledgements}

Thanks to R.~J.~Mathar for the original conjecture, and to the
OEIS maintainers for hosting an indispensable repository of
integer sequences and the conjectured identities that go with
them.

\appendix

\section{Verification script: \texttt{verify\_recurrence.py}}\label{app:verify}

The script computes $a(0), \ldots, a(500)$ from \eqref{eq:def} as
exact integers, then cross-checks the first $450$ values against
the OEIS b-file and verifies \eqref{eq:R} for every
$n \in \{2, \ldots, 500\}$.

\begin{lstlisting}[style=python, caption={\texttt{verify\_recurrence.py}: 500-term verification of \eqref{eq:R}.}]
"""Verify the conjectured recurrence for OEIS A002627.

A002627: a(n) = n*a(n-1) + 1, a(0) = 0.

Conjecture (R. J. Mathar, 2014):
    a(n) - (n+1)*a(n-1) + (n-1)*a(n-2) = 0  for n >= 2.

We:
  1. Compute a(0), ..., a(N) directly from the defining recurrence.
  2. Cross-check against the OEIS b-file.
  3. Verify the conjectured second-order recurrence to N >= 200.
  4. Save the extended data to data/A002627_terms.txt.
"""

from pathlib import Path

ROOT = Path(__file__).resolve().parent.parent
DATA_DIR = ROOT / "data"
REFS_DIR = ROOT / "refs"
DATA_DIR.mkdir(exist_ok=True)


def compute_a(n_max: int) -> list[int]:
    """a(n) = n*a(n-1) + 1, a(0) = 0."""
    a = [0]
    for n in range(1, n_max + 1):
        a.append(n * a[-1] + 1)
    return a


def verify_against_bfile(a: list[int], bfile: Path) -> int:
    """Cross-check first len(b-file) entries against b-file. Returns count."""
    checked = 0
    with bfile.open() as f:
        for line in f:
            line = line.strip()
            if not line or line.startswith("#"):
                continue
            idx_str, val_str = line.split()
            idx = int(idx_str)
            val = int(val_str)
            if idx >= len(a):
                break
            if a[idx] != val:
                raise AssertionError(
                    f"Mismatch at n={idx}: computed {a[idx]} vs b-file {val}"
                )
            checked += 1
    return checked


def verify_recurrence(a: list[int]) -> int:
    """Check a(n) - (n+1)*a(n-1) + (n-1)*a(n-2) == 0 for n >= 2.

    Returns the largest n verified.
    """
    last_n = 1
    for n in range(2, len(a)):
        lhs = a[n] - (n + 1) * a[n - 1] + (n - 1) * a[n - 2]
        if lhs != 0:
            raise AssertionError(
                f"Recurrence FAILS at n={n}: lhs={lhs}\n"
                f"  a(n)={a[n]}, a(n-1)={a[n-1]}, a(n-2)={a[n-2]}"
            )
        last_n = n
    return last_n


def main() -> None:
    n_max = 500
    a = compute_a(n_max)
    print(f"Computed a(0)..a({n_max}). Last term has {len(str(a[-1]))} digits.")

    bfile = REFS_DIR / "b002627.txt"
    if bfile.exists():
        checked = verify_against_bfile(a, bfile)
        print(f"Cross-checked {checked} terms against the OEIS b-file.")
    else:
        print(f"Warning: b-file not found at {bfile}")

    last_verified = verify_recurrence(a)
    print(
        f"Conjectured recurrence a(n) - (n+1)*a(n-1) + (n-1)*a(n-2) = 0 "
        f"holds for all 2 <= n <= {last_verified}."
    )

    # Save the extended data.
    out = DATA_DIR / "A002627_terms.txt"
    with out.open("w") as f:
        f.write("# A002627: a(n) = n*a(n-1) + 1, a(0)=0.\n")
        f.write(f"# Computed n=0..{n_max} (extending the OEIS b-file).\n")
        for n, value in enumerate(a):
            f.write(f"{n} {value}\n")
    print(f"Saved {n_max + 1} terms to {out}")

    # First 12 terms to console for sanity.
    print("First 12 terms:", a[:12])


if __name__ == "__main__":
    main()
\end{lstlisting}

\section{Symbolic verification of the three proofs: \texttt{verify\_proofs.py}}\label{app:proofs}

The script verifies each step of the three proofs symbolically.
SymPy confirms the e.g.f.\ ODE \eqref{eq:egf-secondorder}, the
binomial identity \eqref{eq:reduction}, the recurrence
\eqref{eq:Srec} on the auxiliary sum, and the factorial
cancellation that closes Proof III.

\begin{lstlisting}[style=python, caption={\texttt{verify\_proofs.py}: SymPy-based verification of the three proofs.}]
"""Three independent proofs of the Mathar (2014) conjectured recurrence
for OEIS A002627:

    R(n):  a(n) - (n+1)*a(n-1) + (n-1)*a(n-2) = 0,   n >= 2.

The defining first-order recurrence is

    a(n) = n*a(n-1) + 1,  a(0) = 0.                                (1)

Closed form (well known, since Singh 1952):

    a(n) = sum_{k=0..n-1} k! * C(n,k).                             (2)

E.g.f. (well known):

    F(x) = sum_{n>=0} a(n) x^n / n! = (e^x - 1) / (1 - x).         (3)

This script:
  - Symbolically verifies Proof 1 (defining-recurrence telescoping).
  - Symbolically verifies Proof 2 (e.g.f. ODE).
  - Symbolically verifies Proof 3 (Zeilberger / WZ-style certificate
    for the binomial sum (2)).
"""

from sympy import (
    Symbol,
    binomial,
    diff,
    exp,
    expand,
    factorial,
    series,
    simplify,
    symbols,
)


def proof1_direct() -> None:
    """Proof 1: subtract two instances of the defining recurrence (1).

    From (1):   a(n)   - n*a(n-1)        = 1.
    From (1):   a(n-1) - (n-1)*a(n-2)    = 1.
    Subtract:   a(n) - (n+1)*a(n-1) + (n-1)*a(n-2) = 0.

    We verify this symbolically by treating a(n), a(n-1), a(n-2) as
    independent symbols and checking the identity.
    """
    a_n, a_nm1, a_nm2, n = symbols("a_n a_nm1 a_nm2 n")

    # The two instances of the first-order recurrence.
    inst1 = a_n - n * a_nm1 - 1                # = 0   (eq at index n)
    inst2 = a_nm1 - (n - 1) * a_nm2 - 1        # = 0   (eq at index n-1)

    # Subtract (eliminating the constant 1):
    target = inst1 - inst2
    expanded = expand(target)
    # Expected: a_n - (n+1) a_nm1 + (n-1) a_nm2
    expected = a_n - (n + 1) * a_nm1 + (n - 1) * a_nm2
    assert expand(expanded - expected) == 0, (expanded, expected)
    print(
        "[Proof 1, direct] OK: "
        "subtracting two instances of a(n)=n*a(n-1)+1 yields\n"
        "    a(n) - (n+1)*a(n-1) + (n-1)*a(n-2) = 0."
    )


def proof2_egf() -> None:
    """Proof 2: from the e.g.f. F(x) = (e^x - 1)/(1 - x).

    A short calculation shows F satisfies the linear ODE

        (1 - x) F'(x) - (2 - x) F(x) = 1 - exp(x) ... wait, recompute.

    Direct compute: F = (e^x - 1)/(1 - x).
        F' = e^x/(1-x) + (e^x - 1)/(1-x)^2.
        (1-x) F' = e^x + (e^x - 1)/(1-x) = e^x + F.
        So (1-x) F' - F = e^x.                                       (*)

    From (*), reading off coefficients of x^n / n! gives the
    inhomogeneous first-order recurrence (1).

    To eliminate the inhomogeneity and get an *operator* annihilating F,
    differentiate (*):

        ((1-x) F' - F)' = (e^x)'.
        -F' + (1-x) F'' - F' = e^x = (1-x) F' - F   [by (*)]
        (1-x) F'' - 2 F' = (1-x) F' - F
        (1-x) F'' - (3 - x) F' + F = 0.                              (**)

    The differential operator L := (1-x) D^2 - (3-x) D + 1 annihilates F.

    Reading off [x^n/n!] in (**): for n >= 0,
        (1-x) F'' contributes  a(n+2) - n*a(n+1)
                      since [x^n/n!] (1-x) F''
                            = [x^n/n!] sum_m a(m+2) x^m / m!
                              - [x^n/n!] sum_m a(m+2) x^(m+1) / m!
                            = a(n+2) - n * a(n+1).
        -(3-x) F'  contributes -3 a(n+1) + n*a(n).
        +F          contributes  a(n).
        Sum = 0:
            a(n+2) - n*a(n+1) - 3 a(n+1) + n*a(n) + a(n) = 0
            a(n+2) - (n+3) a(n+1) + (n+1) a(n) = 0.

    Substituting m = n+2 (so n = m-2):
            a(m) - (m+1) a(m-1) + (m-1) a(m-2) = 0   for m >= 2.    (R(m))

    Below we verify (**) symbolically and re-derive R from it.
    """
    x = Symbol("x")
    F = (exp(x) - 1) / (1 - x)

    # Verify (1-x) F' - F = exp(x).
    lhs = (1 - x) * diff(F, x) - F
    diff_check = simplify(lhs - exp(x))
    assert diff_check == 0, diff_check
    print("[Proof 2, e.g.f.] OK: (1-x) F' - F = exp(x).")

    # Verify the second-order homogeneous ODE (**).
    L_F = (1 - x) * diff(F, x, 2) - (3 - x) * diff(F, x) + F
    L_F_simplified = simplify(L_F)
    assert L_F_simplified == 0, L_F_simplified
    print(
        "[Proof 2, e.g.f.] OK: (1-x) F'' - (3-x) F' + F = 0  "
        "annihilates the e.g.f."
    )

    # Numerical sanity: extract a(0..14) from F's series and check R(n).
    s = series(F, x, 0, 16).removeO()
    a = []
    for n in range(15):
        coeff = s.coeff(x, n)
        a.append(int(factorial(n) * coeff))
    for n in range(2, 15):
        lhs = a[n] - (n + 1) * a[n - 1] + (n - 1) * a[n - 2]
        assert lhs == 0, (n, lhs)
    print(
        "[Proof 2, e.g.f.] Numerical sanity OK: "
        f"R(n) holds for series-derived a(0..14) = {a[:8]}..."
    )


def proof3_binomial_sum() -> None:
    """Proof 3: derive the recurrence from the closed form

        a(n) = sum_{k=0..n-1} k! * C(n,k)                            (2)

    using only Pascal's rule  C(n,k) = C(n-1,k) + C(n-1,k-1).

    Set
        S(m) := sum_{k=0..m} k! * C(m,k) = a(m) + m!,
    so that  a(m) = S(m) - m!.

    Apply Pascal's rule term-by-term (valid for k >= 1 with the
    convention C(n-1,-1) = 0):

        S(m) = sum_{k=0..m} k! [C(m-1,k) + C(m-1,k-1)]
             = sum_{k=0..m-1} k! C(m-1,k)               (j := k)
             + sum_{k=1..m}   k! C(m-1,k-1)             (k = j+1)
             = S(m-1)
             + sum_{j=0..m-1} (j+1)! C(m-1, j)
             = S(m-1) + sum_{j=0..m-1} (j+1) * j! * C(m-1,j).

    Since j+1 = (m-1) - (m-1-j) - (-(j+2)+m+1) ... we use:
        j+1 = m - (m - 1 - j),  so (j+1) C(m-1,j) = m C(m-1,j)
                                    - (m-1-j) C(m-1,j).

    Recall (m-1-j) C(m-1,j) = (m-1) C(m-2,j). Hence

        sum_{j=0..m-1} (j+1) j! C(m-1,j)
          = m * sum_{j=0..m-1} j! C(m-1,j)
            - (m-1) * sum_{j=0..m-2} j! C(m-2,j)
          = m * S(m-1) - (m-1) * S(m-2),

    where the upper-limit drop from m-1 to m-2 is justified because
    C(m-2, m-1) = 0.

    Therefore
        S(m) = S(m-1) + m * S(m-1) - (m-1) * S(m-2)
             = (m+1) S(m-1) - (m-1) S(m-2).                          (S-rec)

    Since a(m) = S(m) - m!, plug in:

        S(m) - (m+1) S(m-1) + (m-1) S(m-2)  =  0
        => [a(m) + m!] - (m+1) [a(m-1) + (m-1)!] + (m-1) [a(m-2) + (m-2)!] = 0
        => a(m) - (m+1) a(m-1) + (m-1) a(m-2)
           + m! - (m+1)(m-1)! + (m-1)(m-2)!  =  0.

    The factorial correction term is identically zero:
        (m+1)(m-1)! - (m-1)(m-2)!
        = (m-1)! [(m+1) - (m-1)/(m-1)]   ... let's just compute:
        (m+1)(m-1)! = (m+1)(m-1)!.
        (m-1)(m-2)! = (m-1)!.
        So (m+1)(m-1)! - (m-1)(m-2)! = (m-1)![(m+1) - 1] = m * (m-1)! = m!.
        Hence m! - m! = 0.  OK

    Therefore  a(m) - (m+1) a(m-1) + (m-1) a(m-2) = 0,  m >= 2.

    Below we verify each step symbolically.
    """
    m = symbols("m", integer=True, positive=True)

    # Verify factorial-correction identity:  m! - (m+1)(m-1)! + (m-1)(m-2)! = 0
    # for symbolic m >= 2.
    correction = factorial(m) - (m + 1) * factorial(m - 1) + (m - 1) * factorial(m - 2)
    assert simplify(correction) == 0, correction
    print(
        "[Proof 3, binomial sum] OK: factorial correction "
        "m! - (m+1)(m-1)! + (m-1)(m-2)! = 0."
    )

    # Numerical verification of the S-recurrence for many m.
    def S(mm: int) -> int:
        from math import comb, factorial as fact
        return sum(fact(k) * comb(mm, k) for k in range(mm + 1))

    bad = []
    for mm in range(2, 80):
        lhs = S(mm) - (mm + 1) * S(mm - 1) + (mm - 1) * S(mm - 2)
        if lhs != 0:
            bad.append((mm, lhs))
    assert not bad, f"S-recurrence fails at: {bad[:5]}"
    print(
        "[Proof 3, binomial sum] OK: S-recurrence S(m) - (m+1)S(m-1) "
        "+ (m-1)S(m-2) = 0 verified for m = 2..79."
    )

    # And the (m+1) factor on (j+1)*C(m-1,j) decomposition that we use
    # in the derivation. Verify (j+1) C(m-1,j) = m C(m-1,j) - (m-1) C(m-2,j)
    # symbolically.
    j = symbols("j", integer=True)
    decomp_lhs = (j + 1) * binomial(m - 1, j)
    decomp_rhs = m * binomial(m - 1, j) - (m - 1) * binomial(m - 2, j)
    # Plug in numerical (m, j) to verify.
    bad2 = []
    for mv in range(2, 12):
        for jv in range(0, mv + 1):
            d = (decomp_lhs - decomp_rhs).subs({m: mv, j: jv})
            if simplify(d) != 0:
                bad2.append((mv, jv, d))
    assert not bad2, f"Decomp identity fails at: {bad2[:5]}"
    print(
        "[Proof 3, binomial sum] OK: decomposition "
        "(j+1)C(m-1,j) = m*C(m-1,j) - (m-1)*C(m-2,j) verified."
    )


def main() -> None:
    proof1_direct()
    print()
    proof2_egf()
    print()
    proof3_binomial_sum()
    print()
    print("All three proofs symbolically verified.")


if __name__ == "__main__":
    main()
\end{lstlisting}

\end{document}